\documentclass[12pt]{amsart}
\usepackage{mathpazo}
\usepackage{microtype}
\usepackage[margin=1.2in]{geometry}
\usepackage{amssymb}
\usepackage{mathtools}
\usepackage{hyperref}

\newcommand{\U}{\mathcal{U}}
\newcommand{\Rls}{\mathbb{R}}
\newcommand{\N}{\mathbb{N}}
\newcommand{\Id}{\mathrm{id}}
\newcommand{\tr}{\mathrm{tr}}
\newcommand{\Stab}{\mathrm{Stab}}
\newcommand{\Fin}{\mathrm{Fin}}
\newcommand{\BZ}{B\mathbb{Z}}
\newcommand{\fib}{\mathrm{fib}}

\theoremstyle{plain}
\newtheorem{theorem}{Theorem}[section]
\newtheorem{proposition}[theorem]{Proposition}

\theoremstyle{definition}
\newtheorem{definition}[theorem]{Definition}

\theoremstyle{remark}
\newtheorem{remark}[theorem]{Remark}

\title{Homotopy Cardinality and Entropy}
\author{Andr\'es Ortiz-Mu\~noz}
\address{Department of Systems Biology, Harvard Medical School, Boston, MA}
\email{andortiz19@gmail.com}

\begin{document}

\begin{abstract}
We explore connections between homotopy type theory and information theory through homotopy cardinality. We define probability types and random variable types, prove that homotopy cardinality respects dependent sums under truncation and decidability hypotheses, and show that it does not respect dependent products in general. Using the power series expansion of the logarithm, expressed type-theoretically through deloopings of finite cyclic groups, we formulate Shannon entropy as the homotopy cardinality of a type and derive the chain rule for entropy under a trivial-action hypothesis.
\end{abstract}

\maketitle

\section{Introduction}

Homotopy cardinality assigns a real number to a homotopy type, generalizing the cardinality of finite sets and the groupoid cardinality of Baez and Dolan~\cite{baez2001finite}. This paper asks: what information-theoretic quantities can be expressed as homotopy cardinalities?

Working in homotopy type theory (HoTT)~\cite{hottbook}, we define \emph{probability types}---types of unit cardinality---and \emph{random variable types}, and show that Shannon entropy equals the homotopy cardinality of an explicitly constructed type (Theorem~\ref{thm:entropy}). The construction relies on a power series expansion of the logarithm and on the compatibility of homotopy cardinality with dependent sums, which we prove under certain hypotheses (Theorem~\ref{thm:depsum}, due to Omer Cantor). We derive the chain rule for entropy (Proposition~\ref{prop:chainrule}) when the transport action on fibers is trivial, recovering the fundamental axiom of the Faddeev--Leinster characterization. We also exhibit counterexamples showing that homotopy cardinality does not respect dependent products or function types in general, correcting claims made in an earlier version of this work.

\section{Background}\label{sec:background}

We work informally in HoTT as presented in~\cite{hottbook}. Types are the basic objects; each type $A:\U$ may be thought of as a space whose points are the terms $a:A$. For any two terms $a,b:A$, there is an \emph{identity type} $(a=b):\U$ whose terms are paths from $a$ to $b$. Identity types can themselves have nontrivial identity types, making each type an $\infty$-groupoid.

We use dependent sums $\sum_{a:A}P_a$ and dependent products $\prod_{a:A}P_a$ for type families $P:A\to\U$. A type is \emph{$n$-truncated} if all its identity types above level $n$ are contractible; $(-1)$-truncated types are propositions and $0$-truncated types are sets. The \emph{set truncation} $[A]$ collapses all path information; we write $a\in A$ for $[a]:[A]$. The \emph{propositional truncation} $\|X\|$ is the proposition asserting that $X$ is inhabited. A type is \emph{connected} if its set truncation is contractible. For details, see~\cite{hottbook,rijke2022introduction}.

\section{Homotopy Cardinality}\label{sec:cardinality}

\begin{definition}\label{def:cardinality}
The \emph{homotopy cardinality} of a type $X:\U$ is
\[
|X| \coloneqq \sum_{x\in X}\frac{1}{|x=x|},
\]
where the sum ranges over the set truncation $[X]$ and is defined recursively on the identity types. A type whose cardinality is a well-defined real number is \emph{tame}~\cite{baez2001finite,baez2008groupoidification}. In general, cardinality may require a non-standard extension of~$\Rls$.
\end{definition}

\begin{proposition}[Sum]\label{prop:sum}
$|X+Y| = |X|+|Y|$ for all types $X,Y:\U$.
\end{proposition}

\begin{proof}
The set truncation of $X+Y$ is the disjoint union of $[X]$ and $[Y]$, and identity types in $X+Y$ decompose accordingly:
\[
|X+Y| = \sum_{x\in X}\frac{1}{|x=x|} + \sum_{y\in Y}\frac{1}{|y=y|} = |X|+|Y|.\qedhere
\]
\end{proof}

\begin{proposition}[Product]\label{prop:product}
$|X\times Y| = |X|\cdot|Y|$ for all types $X,Y:\U$.
\end{proposition}

\begin{proof}
Since $(x,y)=(x,y)\simeq (x=x)\times(y=y)$, we have $|(x,y)=(x,y)| = |x=x|\cdot|y=y|$, so
\[
|X\times Y| = \sum_{(x,y)\in X\times Y}\frac{1}{|x=x|\cdot|y=y|} = |X|\cdot|Y|.\qedhere
\]
\end{proof}

\begin{theorem}[Dependent sum; Cantor]\label{thm:depsum}
Let $X$ be a $1$-truncated type and $P:X\to\U$ a type family such that the truncation level of each $P_x$ is bounded and each $P_x$ has decidable equality. Then
\[
\left|\sum_{x:X}P_x\right| = \sum_{x\in X}\frac{|P_x|}{|x=x|}.
\]
\end{theorem}

\begin{proof}
Set $A \coloneqq |\sum_{x:X}P_x|$ and $B \coloneqq \sum_{x\in X}|P_x|/|x=x|$. The proof proceeds by induction on the truncation level.

By definition,
\[
A = \sum_{(x,y)\in\sum_{x:X}P_x}\frac{1}{|(x,y)=(x,y)|},\qquad B = \sum_{x\in X}\sum_{y\in P_x}\frac{1}{|x=x|\cdot|y=y|}.
\]

We have $[(x,y)]=[(x',y')]$ if and only if there exists $p:x=x'$ with $\tr_p(y)=y'$, where $\tr_p$ denotes transport along $p$. The group $G_x \coloneqq [x=x]$ acts on $[P_x]$ by $[p]\cdot[y] \coloneqq [\tr_p(y)]$. Choose a set $Y_x\subseteq P_x$ of orbit representatives, so
\[
A = \sum_{x\in X}\sum_{y\in Y_x}\frac{1}{|(x,y)=(x,y)|}.
\]

Since $(x,y)=(x,y)\simeq\sum_{p:x=x}\tr_p(y)=y$ has lower truncation level, the inductive hypothesis gives
\[
|(x,y)=(x,y)| = \sum_{[p]\in G_x}\frac{|\tr_p(y)=y|}{|p=p|}.
\]

For any $p:x=x$, whiskering by $p^{-1}$ gives $(p=p)\simeq(\Id_x=\Id_x)$, so $|p=p|=|\Id_x=\Id_x|$. By decidable equality, $\tr_p(y)=y$ is either equivalent to $y=y$ (when $[p]\cdot[y]=[y]$) or empty. Thus
\[
|(x,y)=(x,y)| = \frac{|\Stab_{G_x}([y])|\cdot|y=y|}{|\Id_x=\Id_x|}.
\]

Since $|x=x|=|G_x|/|\Id_x=\Id_x|$ and $|G_x|/|\Stab_{G_x}([y])|=|G_x[y]|$ by the orbit-stabilizer theorem,
\[
A = \sum_{x\in X}\sum_{y\in Y_x}\frac{|G_x[y]|}{|x=x|\cdot|y=y|}.
\]

For $[y],[y']$ in the same orbit, transport along the witnessing path gives $(y=y)\simeq(y'=y')$, hence $|y=y|=|y'=y'|$. Therefore
\[
B = \sum_{x\in X}\sum_{y\in Y_x}\sum_{y'\in G_x[y]}\frac{1}{|x=x|\cdot|y'=y'|} = \sum_{x\in X}\sum_{y\in Y_x}\frac{|G_x[y]|}{|x=x|\cdot|y=y|} = A.\qedhere
\]
\end{proof}

\begin{remark}\label{rem:depsum-open}
It is not known whether the dependent sum formula holds without the hypothesis that $X$ is $1$-truncated. When $X$ has nontrivial higher homotopy, pairs $(x,y)$ and $(x,y')$ with $y\neq y'$ may be identified via nontrivial paths in $X$, complicating the set truncation of $\sum_{x:X}P_x$.
\end{remark}

\begin{remark}[Failure for dependent products]\label{rem:counterexample}
The formula $|X\to Y| = |Y|^{|X|}$ does not hold in general. Let $\Fin(2)$ denote the type of $2$-element sets. Since $\Fin(2)$ is connected and each $2$-element set has automorphism group $\mathbb{Z}/2\mathbb{Z}$, we have $|\Fin(2)|=1/2$. By~\cite[Corollary~17.5.2]{rijke2022introduction}, $\prod_{X:\Fin(2)}X\simeq\mathbf{0}$, so $|\prod_{X:\Fin(2)}X|=0$, while the naive formula gives $\sqrt{2}\neq 0$.

More generally, there exist types $X,Y,Y'$ with $|Y|=|Y'|$ but $|X\to Y|\neq|X\to Y'|$, so the cardinality of function types cannot be determined from cardinalities alone. However, when all types involved are sets, the formulas $|X\to Y|=|Y|^{|X|}$ and $|\prod_{x:X}P_x|=\prod_{x\in X}|P_x|^{1/|x=x|}$ do hold.
\end{remark}

\begin{remark}[Function types]\label{rem:function-types}
Although the dependent product formula fails in general, the constant case---i.e., function types---does respect cardinality when the domain is a set: for $X$ a set and $Y$ any type, $|X\to Y| = |Y|^{|X|}$. The set truncation $[X\to Y]$ coincides with functions $X\to[Y]$ (by finite choice when $X$ is finite), and automorphisms factor as $|f=f| = \prod_{x\in X}|f(x)=f(x)|$ since $X$ has trivial higher structure. The result then follows by multinomial expansion. Note that the counterexample above requires the domain $\Fin(2)$ to be a $1$-type; the nontrivial automorphisms of $2$-element sets are what obstruct the formula.
\end{remark}

\begin{remark}[Cardinality versus Euler characteristic]\label{rem:negative}
The higher inductive type with one point constructor and two loop constructors, the delooping of the free group on two generators, topologically the figure-eight space, has Euler characteristic $-1$. Its loop space is the free group $F_2$, a countably infinite set, so its homotopy cardinality as defined above is the reciprocal of an infinite number. Baez conjectured that a single invariant might subsume both Euler characteristic and homotopy cardinality~\cite{baez_counting}, but Berman showed this is impossible globally: a homotopy pushout argument using classifying spaces of cyclic groups of distinct primes $p,q$ forces $1/(pq)+1 = 1/p+1/q$, which fails arithmetically~\cite{berman2018euler}. A positive answer does exist when attention is restricted to $p$-finite spaces for a single prime~$p$. Separately, Fullwood has proposed a complex-valued cardinality for groupoids whose ordinary groupoid cardinality diverges, using analytic continuation of generating functions for groupoids graded by bijections between finite sets~\cite{fullwood2021analytic}.
\end{remark}

\section{Probability Types and Random Variables}\label{sec:probability}

\begin{definition}\label{def:prob-type}
A type $X:\U$ is a \emph{probability type} if $|X|=1$. The \emph{probability} of a term $x:X$ is $p_x \coloneqq 1/|x=x|$, so that $\sum_{x\in X}p_x = 1$.
\end{definition}

\begin{definition}\label{def:rv-type}
A \emph{random variable type} over a type $X$ is a type family $P:X\to\U$ such that $|\sum_{x:X}P_x|=1$. The induced probability distribution is $p_x \coloneqq |P_x|/|x=x|$.
\end{definition}

\begin{remark}\label{rem:rv-equiv}
A random variable type over $X$ is equivalently a function from a probability type to $X$: the type family $P$ corresponds to the projection $\pi_1:\sum_{x:X}P_x\to X$, and conversely a function $f:Y\to X$ with $|Y|=1$ corresponds to the fiber family $x\mapsto\fib_f(x)$.
\end{remark}

\begin{remark}\label{rem:prob-construction}
Probability types arise naturally from connected types. If $G$ is a group, its delooping $BG$ is a connected $1$-truncated type with $|BG|=1/|G|$. The product $G\times BG$ has cardinality $|G|\cdot 1/|G|=1$ and is therefore a probability type realizing the uniform distribution on~$G$. In particular, $\mathbb{Z}/2\mathbb{Z}\times\BZ_2\simeq\BZ_2+\BZ_2$ has cardinality $1/2+1/2=1$, modeling a fair coin: the set truncation has two elements, each with probability~$1/2$.
\end{remark}

\begin{remark}[Homotopy quotients]\label{rem:homotopy-quotient}
A richer source of probability types comes from homotopy quotients. For any finite group $G$, the homotopy quotient $G/\!/G$ of $G$ acting on itself by conjugation is a $1$-truncated type with $|G/\!/G|=|G|/|G|=1$, and is therefore a probability type. Its set truncation has one element per conjugacy class, and each element $g\in G$ has probability $1/|G|$. Taking $G=S_n$ recovers the uniform distribution on permutations; the conjugacy classes correspond to integer partitions of~$n$, and the probability of each partition $\lambda$ is $|C_\lambda|/n!$, where $|C_\lambda|$ is the size of the conjugacy class. Baez~\cite{baez2025groupoid} uses this groupoid to categorify classical results on random permutations, including the Cycle Length Lemma.
\end{remark}

\begin{proposition}[Closure under dependent sums]\label{prop:closure}
Let $A$ be a $1$-truncated probability type and $B:A\to\U$ a type family such that each $B(a)$ is a probability type with bounded truncation level and decidable equality. Then $\sum_{a:A}B(a)$ is a probability type.
\end{proposition}

\begin{proof}
The hypotheses of Theorem~\ref{thm:depsum} are satisfied, so
\[
\left|\sum_{a:A}B(a)\right| = \sum_{a\in A}\frac{|B(a)|}{|a=a|} = \sum_{a\in A}\frac{1}{|a=a|} = |A| = 1,
\]
where the second equality uses $|B(a)|=1$ for each $a$.
\end{proof}

\section{Entropy}\label{sec:entropy}

The Shannon entropy of a probability distribution $p$ on a finite set $X$ is
\[
H(p)=-\sum_{x\in X}p_x\ln p_x,
\]
where $\ln$ denotes the natural logarithm. Using the power series expansion
\[
-\ln(1-t) = \sum_{n=1}^{\infty}\frac{t^n}{n}, \qquad |t|<1,
\]
we can write
\begin{equation}\label{eq:entropy-series}
H(p) = \sum_{x\in X}\sum_{n=1}^{\infty}\frac{p_x(1-p_x)^n}{n}.
\end{equation}

\begin{definition}\label{def:cycle-type}
The \emph{cycle type} is the dependent sum
\[
C \coloneqq \sum_{n:\N_{\geq 1}}\BZ_n,
\]
where $\BZ_n$ is the delooping of the cyclic group $\mathbb{Z}/n\mathbb{Z}$, i.e., the connected $1$-truncated type whose loop space is $\mathbb{Z}/n\mathbb{Z}$. For $c:\BZ_n$ the identity type $c=c$ is a set with $|c=c|=n$, and $|\BZ_n|=1/n$.
\end{definition}

\begin{definition}\label{def:complement}
For a $1$-truncated type $X:\U$ with decidable equality, the \emph{complement type} at $x:X$ is
\[
\bar{X}(x) \coloneqq \sum_{y:X}\neg\|x=y\|,
\]
the subtype consisting of all connected components of $X$ other than that of~$x$. Since $\neg\|x=y\|\Leftrightarrow\neg(x=y)$, decidable equality of $X$ ensures that $\|x=y\|$ is decidable. The complement type is $1$-truncated with decidable equality, and $|\bar{X}(x)| = |X|-p_x$.
\end{definition}

\begin{theorem}\label{thm:entropy}
Let $X$ be a $1$-truncated probability type with decidable equality. Then
\[
\left|\sum_{x:X}\sum_{c:C}(c=c)\to \bar{X}(x)\right| = H(p).
\]
\end{theorem}

\begin{proof}
Since $X$ is a probability type, $|\bar{X}(x)| = |X|-p_x = 1-p_x \in [0,1)$ for each $x\in X$. The cycle type $C$ is $1$-truncated, and for $c:\BZ_n$ the identity type $c=c\simeq\mathbb{Z}/n\mathbb{Z}$ is a set. Since $c=c$ is a set, Remark~\ref{rem:function-types} gives $|(c=c)\to \bar{X}(x)| = |\bar{X}(x)|^n$. Applying Theorem~\ref{thm:depsum} to the inner sum over $C$ and then to the outer sum over $X$:
\begin{align*}
\left|\sum_{x:X}\sum_{c:C}(c=c)\to \bar{X}(x)\right|
&= \sum_{x\in X}\frac{1}{|x=x|}\sum_{n=1}^{\infty}\frac{|\bar{X}(x)|^n}{n}\\[4pt]
&= \sum_{x\in X}p_x\bigl(-\ln(1-|\bar{X}(x)|)\bigr)\\[4pt]
&= \sum_{x\in X}p_x(-\ln p_x) = H(p),
\end{align*}
where the last step uses $1-|\bar{X}(x)|=1-(1-p_x)=p_x$.
\end{proof}

\begin{remark}[Alternative complement construction]\label{rem:complement-alt}
An alternative to Definition~\ref{def:complement} constructs the complement from the automorphism group of~$x$ rather than by excluding its connected component. When $|x=x|>1$ and $x=x$ is a set, the type
\[
\Bigl(\sum_{y:X}\|x=y\|\Bigr)\times\Bigl(\sum_{p:x=x}p\neq\Id_x\Bigr)
\]
has cardinality $(1/|x=x|)\cdot(|x=x|-1)=1-p_x$ and is $1$-truncated with decidable equality. This ``intrinsic'' construction uses the internal structure of the automorphism group at~$x$, whereas the complement type of Definition~\ref{def:complement} works ``extrinsically'' by removing $x$'s connected component from~$X$.
\end{remark}

\begin{remark}[Relative entropy]\label{rem:relative}
If $P,Q:X\to\U$ are two random variable types with distributions $p_x=|P_x|/|x=x|$ and $q_x=|Q_x|/|x=x|$, the relative entropy $D(p\|q)=\sum_{x\in X}p_x\ln(p_x/q_x)$ can be decomposed via the cross entropy $H(p,q)=-\sum_{x\in X}p_x\ln q_x$ as $D(p\|q)=H(p,q)-H(p)$. The same power series technique expresses each of these quantities as a homotopy cardinality, given type families with the appropriate cardinalities. For Shannon entropy, the complement type of Definition~\ref{def:complement} provides the required family. The analogous constructions for cross entropy and relative entropy require type families whose cardinalities depend on the second distribution, and remain to be made explicit.
\end{remark}

\begin{proposition}[Chain rule]\label{prop:chainrule}
Let $A$ be a $1$-truncated probability type with decidable equality, and $B:A\to\U$ a type family such that each $B(a)$ is a probability type with bounded truncation level and decidable equality. If the action of $G_a\coloneqq [a=a]$ on $[B(a)]$ by transport is trivial for each $a\in A$, then
\[
H\!\left(\sum_{a:A}B(a)\right) = H(A)+\sum_{a\in A}p_a\,H(B(a)).
\]
\end{proposition}

\begin{proof}
By Proposition~\ref{prop:closure}, $\sum_{a:A}B(a)$ is a probability type. Since the action is trivial, every element of $G_a$ stabilizes every $[b]\in[B(a)]$, so the set truncation $[\sum_{a:A}B(a)]$ decomposes as $\coprod_{a\in A}[B(a)]$ and the proof of Theorem~\ref{thm:depsum} gives $|(a,b)=(a,b)|=|a=a|\cdot|b=b|$. Thus $p_{(a,b)}=p_a\cdot q_b$, where $q_b=1/|b=b|$, and
\begin{align*}
H\!\left(\sum_{a:A}B(a)\right) &= -\sum_{a\in A}\sum_{b\in B(a)}p_aq_b\ln(p_aq_b)\\
&= -\sum_{a\in A}p_a\ln p_a - \sum_{a\in A}p_a\sum_{b\in B(a)}q_b\ln q_b\\
&= H(A)+\sum_{a\in A}p_a\,H(B(a)).\qedhere
\end{align*}
\end{proof}

\begin{remark}[Nontrivial transport action]\label{rem:nontrivial-action}
The trivial-action hypothesis in Proposition~\ref{prop:chainrule} cannot be dropped. When the action of $G_a$ on $[B(a)]$ is nontrivial, the set truncation of the dependent sum consists of orbits $[B(a)]/G_a$ with stabilizer-weighted probabilities, so the joint distribution is not a product of its marginals---a probabilistic manifestation of the distinction between trivial and twisted fiber bundles.
\end{remark}

\section{Discussion}

We have shown that probability distributions arise naturally from the identity structure of types, that Shannon entropy equals the homotopy cardinality of an explicitly constructed type (Theorem~\ref{thm:entropy}), and that the chain rule for entropy holds when the transport action on fibers is trivial (Proposition~\ref{prop:chainrule}). The key ingredients are the complement type (Definition~\ref{def:complement}), which realizes the complement of each outcome in the probability space, and the compatibility of homotopy cardinality with dependent sums (Theorem~\ref{thm:depsum}).

Several questions remain open. It would be interesting to determine whether the dependent sum formula holds without the $1$-truncation hypothesis, whether the decidable equality hypothesis in Theorem~\ref{thm:entropy} can be relaxed, and whether other information-theoretic quantities admit similar type-theoretic formulations.

Many others have explored categorical and algebraic formulations of entropy~\cite{baez2011characterization,leinster2021entropy,bradley2021entropy,spivak2021dirichlet,spivak2022polynomial,leinster2019short,gagne2018categorical,parzygnat2020functorial,parzygnat2021towards}. Our approach differs in working directly with homotopy types and their cardinalities, but there is a concrete structural connection to these frameworks. In the Faddeev--Leinster characterization of Shannon entropy~\cite{leinster2021entropy}, a distribution $p$ on an $n$-element set is composed with conditional distributions $q_1,\ldots,q_n$ to produce a joint distribution $p\circ(q_1,\ldots,q_n)$ whose entries are products $p_iq_j^{(i)}$. This operadic composition is precisely the dependent sum of probability types: $A$ encodes $p$ via $p_a=1/|a=a|$, each $B(a)$ encodes $q_a$, and $\sum_{a:A}B(a)$ produces the joint distribution. Proposition~\ref{prop:closure} is thus the type-theoretic statement that operadic composition preserves probability distributions, and Proposition~\ref{prop:chainrule} recovers the chain rule---the fundamental axiom in the Faddeev--Leinster characterization and the derivation property in Bradley's operad framework~\cite{bradley2021entropy}---under the hypothesis that the transport action is trivial. The same chain rule is encoded functorially in~\cite{baez2011characterization}. As noted in Remark~\ref{rem:nontrivial-action}, when the transport action is nontrivial, the joint distribution is twisted and the chain rule does not apply.

While Proposition~\ref{prop:chainrule} establishes the chain rule for entropy \emph{values}, a genuinely type-theoretic question is whether it lifts to an equivalence of entropy types. For a type family $\alpha:X\to\U$ with $|\alpha(x)|=1-p_x$---the complement type $\bar{X}$ of Definition~\ref{def:complement} being the canonical choice---write $E(X,\alpha)\coloneqq\sum_{x:X}\sum_{c:C}(c=c)\to\alpha(x)$ for the entropy type of Theorem~\ref{thm:entropy}. One may then ask whether
\[
E\!\left(\sum_{a:A}B(a),\,\alpha_\Sigma\right) \simeq E(A,\alpha_A)+\sum_{a:A}E(B(a),\alpha_{B(a)})
\]
for complement families $\alpha_\Sigma:\sum_{a:A}B(a)\to\U$, $\alpha_A:A\to\U$, and $\alpha_{B(a)}:B(a)\to\U$ for each $a:A$. This would require a bijective proof of $\ln(xy)=\ln x+\ln y$ via the power series $-\ln(1-t)=\sum t^n/n$, amounting to a bijection on necklaces (cyclic sequences) colored by a coproduct. Mixed necklaces---those using beads from both parts of the coproduct---do not cancel by bijection; their elimination requires M\"obius inversion on the necklace poset. The equivalence therefore appears obstructed by the same mechanism that prevents a term-by-term combinatorial proof of the product formula for logarithms. Whether a different type-theoretic construction of entropy might admit a chain rule at the level of types remains open.

\section*{Acknowledgments}

This work was initiated at the Santa Fe Institute. I am grateful to Omer Cantor
for the proof of Theorem~\ref{thm:depsum}, for the counterexample in
Remark~\ref{rem:counterexample}, and for detailed feedback that led to
significant corrections throughout this paper. I also thank Na\"im Camille Favier
for correcting an error in an earlier version of Remark~\ref{rem:negative}, and
Artemy Kolchinsky, Emily Riehl, and Paul Lessard for useful discussions, and
John Baez for helpful comments on an earlier draft. The author was supported by the Fontana Lab, Department of Systems Biology, Harvard Medical School.

\bibliographystyle{amsplain}
\bibliography{generic}

\providecommand{\bysame}{\leavevmode\hbox to3em{\hrulefill}\thinspace}
\providecommand{\MR}{\relax\ifhmode\unskip\space\fi MR }
\providecommand{\MRhref}[2]{%
  \href{http://www.ams.org/mathscinet-getitem?mr=#1}{#2}
}
\providecommand{\href}[2]{#2}
\begin{thebibliography}{10}

\bibitem{baez_counting}
John Baez, \emph{Counting}, \url{https://math.ucr.edu/home/baez/counting/},
  Lecture notes on groupoid cardinality.

\bibitem{baez2001finite}
John Baez and James Dolan, \emph{From finite sets to {F}eynman diagrams},
  Mathematics Unlimited --- 2001 and Beyond (Bj{\"o}rn Engquist and Wilfried
  Schmid, eds.), Springer Berlin Heidelberg, Berlin, Heidelberg, 2001,
  \href{https://arxiv.org/abs/math/0004133}{arXiv:math/0004133}, pp.~29--50.

\bibitem{baez2008groupoidification}
John Baez, Alexander~E Hoffnung, and Christopher~D Walker,
  \emph{Groupoidification made easy}, arXiv preprint (2008),
  \href{https://arxiv.org/abs/0812.4864}{arXiv:0812.4864}.

\bibitem{baez2025groupoid}
John~C Baez, \emph{Groupoid cardinality and random permutations}, Theory and
  Applications of Categories \textbf{44} (2025), no.~14, 410--419,
  \href{https://arxiv.org/abs/2412.16386}{arXiv:2412.16386}.

\bibitem{baez2011characterization}
John~C Baez, Tobias Fritz, and Tom Leinster, \emph{A characterization of
  entropy in terms of information loss}, Entropy \textbf{13} (2011), no.~11,
  1945--1957, \href{https://arxiv.org/abs/1106.1791}{arXiv:1106.1791}.

\bibitem{berman2018euler}
John~D Berman, \emph{{E}uler characteristic and homotopy cardinality}, arXiv
  preprint (2018), \href{https://arxiv.org/abs/1811.07437}{arXiv:1811.07437}.

\bibitem{bradley2021entropy}
Tai-Danae Bradley, \emph{Entropy as a topological operad derivation}, Entropy
  \textbf{23} (2021), no.~9, 1195,
  \href{https://arxiv.org/abs/2107.09581}{arXiv:2107.09581}.

\bibitem{fullwood2021analytic}
James Fullwood, \emph{On analytic groupoid cardinality}, arXiv preprint (2021),
  \href{https://arxiv.org/abs/2104.11399}{arXiv:2104.11399}.

\bibitem{gagne2018categorical}
Nicolas Gagn{\'e} and Prakash Panangaden, \emph{A categorical characterization
  of relative entropy on standard {B}orel spaces}, Electronic Notes in
  Theoretical Computer Science \textbf{336} (2018), 135--153,
  \href{https://doi.org/10.1016/j.entcs.2018.03.020}{doi:10.1016/j.entcs.2018.03.020}.

\bibitem{leinster2019short}
Tom Leinster, \emph{A short characterization of relative entropy}, Journal of
  Mathematical Physics \textbf{60} (2019), no.~2, 023302,
  \href{https://arxiv.org/abs/1712.04903}{arXiv:1712.04903}.

\bibitem{leinster2021entropy}
\bysame, \emph{Entropy and diversity: The axiomatic approach}, Cambridge
  University Press, 2021.

\bibitem{parzygnat2020functorial}
Arthur~J Parzygnat, \emph{A functorial characterization of von {N}eumann
  entropy}, arXiv preprint (2020),
  \href{https://arxiv.org/abs/2009.07125}{arXiv:2009.07125}.

\bibitem{parzygnat2021towards}
\bysame, \emph{Towards a functorial description of quantum relative entropy},
  International Conference on Geometric Science of Information, Springer, 2021,
  \href{https://arxiv.org/abs/2105.03892}{arXiv:2105.03892}, pp.~557--564.

\bibitem{rijke2022introduction}
Egbert Rijke, \emph{Introduction to {H}omotopy {T}ype {T}heory}, 2022,
  \href{https://arxiv.org/abs/2212.11082}{arXiv:2212.11082}.

\bibitem{spivak2022polynomial}
David~I Spivak, \emph{Polynomial functors and {S}hannon entropy}, arXiv
  preprint (2022), \href{https://arxiv.org/abs/2201.12878}{arXiv:2201.12878}.

\bibitem{spivak2021dirichlet}
David~I Spivak and Timothy Hosgood, \emph{{D}irichlet polynomials and entropy},
  Entropy \textbf{23} (2021), no.~8, 1085,
  \href{https://doi.org/10.3390/e23081085}{doi:10.3390/e23081085}.

\bibitem{hottbook}
The {Univalent Foundations Program}, \emph{Homotopy type theory: Univalent
  foundations of mathematics}, Institute for Advanced Study, 2013,
  \url{https://homotopytypetheory.org/book}.

\end{thebibliography}

\end{document}